\documentclass[12pt]{article}
\usepackage[latin1]{inputenc}
\usepackage{amssymb,amsmath,amsfonts}
\newtheorem{theorem}{Theorem}
\newtheorem{proposition}[theorem]{Proposition}
\newtheorem{lemma}[theorem]{Lemma}

\newtheorem{remark}[theorem]{Remark}

\newcommand{\R}{\mathbb{R}}

\newcommand{\Sf}{\mathbb{S}}
\newcommand{\C}{\mathbb{C}}

\newcommand{\spa}{\mbox{span}}

\newcommand{\po}{{\hspace*{-1ex}}{\bf .  }}
\newcommand{\ii}{isometric immersion }
\newcommand{\iis}{isometric immersions }

\newcommand{\nap}{\nabla^{\perp}}

\newcommand{\nab}{\tilde\nabla}

\def\<{{\langle}}
\def\>{{\rangle}}

\def\n{\nabla}
\def\d{\partial}

\def\a{\alpha}
\def\va{\varphi}

\def\bea{\begin{eqnarray*} }
\def\eea{\end{eqnarray*} }
\def\be{\begin{equation} }
\def\ee{\end{equation} }

\def\nap{\nabla^\perp}
\def\proof{\noindent\emph{Proof: }}
\def\qed{\ifhmode\unskip\nobreak\fi\ifmmode\ifinner
\else\hskip5 pt \fi\fi\hbox{\hskip5 pt \vrule width4 pt
height6 pt  depth1.5 pt \hskip 1pt }}

\begin{document}

\title{A class of complete minimal  submanifolds and
their associated families of genuine deformations}
\author{M. Dajczer and Th. Vlachos}
\date{}
\maketitle

\begin{abstract} 
Concerning the problem of classifying complete submanifolds of Euclidean 
space with codimension two admitting genuine isometric deformations, until
now the only known examples with the maximal possible rank four 
are the real Kaehler  minimal submanifolds classified by Dajczer-Gromoll 
\cite{dg3} in parametric form. These submanifolds behave like minimal surfaces, 
namely, if simple connected  either they admit a nontrivial one-parameter 
associated family of isometric deformations or are holomorphic. 

In this paper, we characterize a new class of complete minimal genuinely 
deformable Euclidean submanifolds of rank four but now the structure of their 
second fundamental  and the way it gets modified while deforming is quite more 
involved than in the Kaehler case. 
This can be seen as a strong indication that the above classification problem 
is quite challenging. Being minimal, the submanifolds we introduced are also 
interesting by themselves. In particular, because associated to any complete
holomorphic curve in $\C^N$ there is such a submanifold and, beside, the
manifold in general is not Kaehler.
\end{abstract}

Some of the very basic question in the local and global theory of isometric 
immersions of Riemannian manifolds into Euclidean space remain in good part 
unanswered. For instance,  outside some special cases it is not known which 
is the lowest codimension for which a given Riemannian manifold admits an 
isometric immersion. On one hand, there are several results that assure that 
a submanifold  must be unique, that is, isometrically rigid, when lying in 
its lowest possible codimension. On the other hand,
there are few theorems classifying isometrically deformable submanifolds 
and their deformations.  This is due to the fact that rigidity is a ``generic"
property while being deformable is certainly not, and hence a situation harder
to deal with.

 The exception for the deformation problem is the case of hypersurfaces.
In fact, in the local case the problem was mostly solved by Sbrana \cite{sb} 
and Cartan \cite{ca} about a century ago; see \cite{dft1} for details and 
a modern presentation. A solution to the problem for compact hypersurfaces 
was given by Sacksteder \cite{sa} and by  Dajczer-Gromoll \cite{dg2} in 
the complete case. But solving the deformation problem in codimension two  
turns out to be very challenging even in the more restrictive case of 
complete manifolds. 

In dealing with the isometric deformation problem in higher codimension, it has
to be taken into account that any submanifold of a deformable submanifold has the
isometric deformations induced by the latter.  In order to obtain classifications, 
it is natural to exclude this type of deformations and only study the remaining ones 
that were called genuine deformations in \cite{df2}. 

An \ii $\hat{f}\colon M^n\to\R^{n+p}$ is a  
\emph{genuine deformation} of a given \ii $f\colon M^n\to\R^{n+p}$, $p\geq 2$, 
if there is no  open subset $U\subset M^n$ along which  $f|_U$ and $\hat{f}|_U$ extend 
isometrically. That  $f\colon M^n\to \R^{n+p}$ and  
$\hat f\colon M^n\to\R^{n+p}$ \emph{extend isometrically} means that there is 
an isometric embedding \mbox{$j\colon M^n\hookrightarrow N^{n+q}$}, $1\leq q<p$,  
into a Riemannian manifold $N^{n+q}$ and there are \iis $F\colon N^m\to\R^{n+q}$ and 
$\hat F\colon N^m\to\R^{n+q}$ such that $f=F\circ j$ and $\hat f=\hat F\circ j$, 
i.e., the following diagram commutes: 
\medskip 
 
\begin{picture}(110,84)   
\put(160,31){$M^n$} 
\put(209,31){$N^{n+q}$} 
\put(413,31){$(1)$} 
\put(242,62){$\R^{n+p}$} 
\put(242,0){$\R^{n+p}$} 
\put(195,59){${}_f$} 
\put(195,9){${}_{\hat f}$} 
\put(240,46){${}_F$} 
\put(239,26){${}_{\hat F}$} 
\put(193,40.5){${}_j$} 
\put(228,42){\vector(1,1){16}} 
\put(228,28){\vector(1,-1){16}} 
\put(177,44){\vector(3,1){60}} 
\put(177,26){\vector(3,-1){60}} 
\put(185,34){\vector(1,0){21}} 
\put(185,36){\oval(7,4)[l]}
\end{picture}
\bigskip

The only general result for submanifolds that admit genuine
deformations known at this time is the local result due to Dajczer-Florit \cite{df2}.
In low codimension, they showed that genuine deformations are only possible for certain
class of ruled submanifolds and gave a lower bound for the dimension of the rulings.
In the special case of codimension two, in order to admit genuine 
deformations a submanifold without flat points must have rank $\rho$ at most 
four at any point. By $\rho$ we denote the rank of the Gauss map, that is, 
$\rho=n-\nu$, where $\nu$ stands for the standard index of relative nullity, 
namely, the dimension of the kernel of the second fundamental form. 

In this paper, we are interested in the global problem of genuine deformations of \iis
with codimension two. In fact, we deal with the noncompact case since for compact 
submanifolds the deformation problem was already solved by 
Dajczer-Gromoll \cite{dg4}. We point out that there exist several local results 
on genuine deformations in the special case of submanifold of rank $\rho=2$ 
but these manifolds  are never complete; see \cite{df2}, \cite{dft2} and \cite{fg}.
In particular, there are the minimal ones that were parametrically classified 
in \cite{df1}. They admit a one-parameter associated family of isometric deformations 
whose geometric nature was recently described in \cite{dv}.

At this time, there is only one classification result on deformations for complete
noncompact submanifolds in Euclidean space with codimension two, namely,
the one given in \cite{dg3} of minimal but non-holomorphic isometric immersions 
of Kaehler manifolds. If simply connected such a submanifold admits 
a nontrivial one-parameter associated family of isometric deformations; see \cite{dg1}.
These submanifolds are ruled (i.e., foliated by complete Euclidean spaces) with
rulings of codimension two and have rank $\rho=4$ almost everywhere.
As in the  case of minimal surfaces, the associated family is obtained by
composing its second fundamental form with an orthogonal parallel tensor in the tangent
bundle given in terms of the complex structure of the manifold. The tensor amounts to  a 
rotation of constant angle while keeping the  the normal bundle and the induced 
connection unchanged.  Basically, this is also the situation of the local case discussed 
to in the preceding  paragraph.

In this paper, we  parametrically construct and characterize a new class of complete
minimal ruled submanifolds that also admit a one-parameter associated family of isometric
deformations.  As  before, the rulings have codimension two and the rank is $\rho=4$ 
almost everywhere. Moreover, the deformations  are obtained while keeping unchanged 
the normal bundle and connection.
But now, the second fundamental form of the deformed submanifold relates to the
initial one in a much more complex form, in particular, no orthogonal tensor
is involved.  

It is an interesting question if the above two families of complete ruled minimal
submanifolds exhaust all examples in the same class that admit genuine deformations.
For instance, they may be examples such that the integral leaf exists but it is 
not totally geodesic. Of course, a much more challenging classification problem 
of complete submanifolds of rank four would be to drop one of the conditions, for
instance being minimal or ruled.
In the Kaehler case, it follows from \cite{dg3} that there are a lot more examples 
without complete rulings. From the recent results in \cite{dv2} it follows that 
this is also the situation in our case.

Finally, we observe that some arguments in this paper involve some unexpected 
long but straightforward computations that will be only sketched.

\section{The $1$-isotropic surfaces}

In this section, we discuss some properties of the $1$-isotropic surfaces in Euclidean 
space that are the basic tool for the construction of the minimal submanifolds 
that are the object of this paper.
\vspace{1,5ex}

Let $g\colon L^2\to\R^{n+2}$ denote an isometric immersion of a
two-dimensional oriented Riemannian manifold into Euclidean space.
The $k^{th}$\emph{-normal space} of $g$ at $p\in L^2$ for $k\geq 1$ is 
given by
$$
N^g_k(p)=\spa\{\a_g^{k+1}(X_1,\ldots,X_{k+1}):X_1,\ldots,X_{k+1}\in T_pL\}
$$
where $\a_g^2=\a_g\colon TL\times TL\to N_gL$ is the standard 
second fundamental form with values in the normal bundle and  
$$
\a_g^s\colon TL\times\cdots\times TL\to N_gL,\;\; s\geq 3, 
$$
is the symmetric tensor called
the $s^{th}$\emph{-fundamental form} defined inductively by
$$
\a_g^s(X_1,\ldots,X_s)=\left(\nabla^\perp_{X_s}\ldots
\nabla^\perp_{X_3}\a_g(X_2,X_1)\right)^\perp.
$$
Here  $\nabla^{\perp}$ is the induced connection in the normal bundle $N_gL$
and $(\;\;)^\perp$ means taking the projection onto the normal complement of 
$N^g_1\oplus\cdots\oplus N^g_{s-2}$ in $N_gL$.

Assume further that $g\colon L^2\to\R^{n+2}$ is minimal and substantial. 
The latter means that the codimension cannot be reduced,
in fact, not even locally since minimal surfaces are real analytic.
Then, on an open dense subset of $L^2$ the normal bundle of $g$ splits as
$$
N_gL=N_1^g\oplus N_2^g\oplus\dots\oplus N_m^g,\;\;\; m=[(n-1)/2],
$$
since all higher normal bundles  have rank two except possible the last 
one that has rank one if $n$ is odd; see \cite{Ch}, \cite{df1} or \cite{Sp} 
for details. 
Moreover, the orientation of $L^2$  induces an orientation on each plane 
vector bundle $N_s^g$  given by the ordered pair 
$$
\xi_1^s=\a_g^{s+1}(X,\ldots,X),\;\;\;\xi_2^s=\a_g^{s+1}(JX,\ldots,X)
$$
where $0\neq X\in TL$ and $J$ is the complex structure of $L^2$ 
determined by the metric and orientation.

If  $L^2$ is simply connected, the generalized Weierstrass parametrization implies 
that there exists a one-parameter \emph{associated family} of minimal immersions; 
see \cite{ho}.  An alternative way to see this goes as follows: for each constant 
$\theta\in\Sf^1=[0,\pi)$ consider the orthogonal parallel tensor field 
$$
J_{\theta}=\cos\theta I+\sin\theta J
$$
where $I$ is the identity map.  Then, the symmetric section $\a_g(J_\theta\cdot,\cdot)$ 
of the bundle $\text{Hom}(TL\times TL,N_g L)$ satisfies the Gauss, Codazzi and Ricci 
equations with respect to the normal bundle and normal connection of $g$; see \cite{dg1}. 
Therefore, there exists an isometric minimal immersion  $g_{\theta}\colon L^2\to\R^{n+2}$ 
whose second fundamental form is
$$
\a_{g_{\theta}}(X,Y)=\phi_\theta\a_g(J_{\theta}X,Y)
$$ 
where $\phi_\theta\colon N_gL\to N_{g_{\theta}}L$ is the parallel 
vector bundle isometry that identifies the normal bundles. 
Explicitly,  the immersion is given by the line integral 
$$
g_\theta(x)=\int_{p_0}^x g_*\circ J_\theta
$$
where $p_0$ is any fixed point in $L^2$. In particular, we have that
$g_{\theta *}=g_*\circ J_\theta$. Thus $\phi_{\theta}$ is nothing else
than parallel identification in $\R^{n+2}$ that identifies all normal subbundles
$N_j^g$ with $N_j^{g_\theta}$, $j\geq 1$, and for simplicity will be dropped from now on.
It turns out that the associated family is \emph{trivial} (i.e.,  each $g_{\theta}$ 
is congruent to $g$) if and only if $g$ is a holomorphic curve with respect to some 
complex structure of the ambient space; cf.\ \cite{df1}. 

\begin{remark}\po {\em The case when $L^2$ above is non-simply-connected was 
considered in \cite{dv2}}.
\end{remark}

Now assume that $g\colon L^2\to\R^{n+2}$, $n\geq 2$, is substantial and 
\emph{$1$-isotropic}.  The latter means that the surface is minimal and 
that the  ellipse of curvature  at all points is 
a circle.  Recall that the \textit{ellipse of curvature}
$\mathcal{E}^g(p)\subset N^g_1(p)$ of $g$ at $p\in L^2$ is defined as
$$
\mathcal{E}^g(p)=\{\a_g(X_{\psi},X_{\psi}) :
X_{\psi}=\cos\psi X+\sin\psi JX\;\;\mbox{and}\;\;\psi\in [0,2\pi)\}
$$
where $X\in T_pL$ has unit length. 
\vspace{1,5ex}

The argument for the following result is basically due to Chern \cite{cr}.

\begin{proposition}\po\label{extend} Let $L_0$ be the open subset of 
$L^2$ where $\dim N_{1}^g(p)=2$.
Then, $L^2\smallsetminus L_0$ consists at most of isolated points and the vector 
bundle $N_1^g|_{L_0}$ extends smoothly to a  plane bundle over $L^2$ still 
denoted by $N^g_1$. 
\end{proposition}

\proof 
The complexified tangent bundle $TL\otimes\C$ decomposes into
the eigenspaces of the complex structure $J$ corresponding to the eigenvalues
$i$ and $-i$ denoted by $T^{\prime}L$ and $T^{\prime\prime}L$, respectively.
The second fundamental form can be complex linearly extended to
$TL\otimes\C$ with values in the complexified vector bundle
$N_gL\otimes\C$ and then decomposed
into its $(p,q)$-components, $p+q=2,$ which are tensor products of $p$ many
1-forms vanishing on $T^{\prime\prime}L$ and $q$ many 1-forms vanishing
on $T^{\prime}L$. Since the surface is minimal the $(1,1)$-part of $\a_g$
vanishes, i.e., $\a_g(\d_z,\bar\d_z)=0$ where $z$ is a complex coordinate.
We thus have the splitting
\be\label{deco}
\a_g= \a^{(2,0)}+ \a^{(0,2)}\;\;\mbox{where}\;\;\a^{(0,2)}=\overline{\a^{(2,0)}}.
\ee
The Codazzi equation implies that
$$
\nap_{\bar\d_z}\a_g(\d_z,\d_z)=0
$$
which means that $\a^{(2,0)}$ is holomorphic as a $ N_gL\otimes\C$-valued
tensor field. 

Since $g$ is $1$-isotropic, then $\dim N_1^g(p_0)<2$ if and only if
$\a_g(p_0)=0$. Moreover, in (\ref{deco})
the summands are perpendicular with respect to the hermitian inner product. 
Hence, the zeros of $\a_g$ are precisely the zeros of $\a^{(2,0)}$. 
Since $\a^{(2,0)}$ is holomorphic, we conclude that its zeros are isolated, 
and hence $L^2\smallsetminus L_0$ consists at most of isolated points.  

Let $(U,z)$ be a complex chart around a point 
$p_0\in L^2\smallsetminus L_0$ with $z(p_0)=0$.
Since $\a^{(2,0)}$ it is not identically zero and $p_0$ is a zero of it,
around $p_0$ we may write
$$
\a^{(2,0)}=z^m \a^{*(2,0)}
$$
for a positive integer $m$, where $\a^{*(2,0)}$ is a tensor field of
type $(2,0)$ with $\a^{*(2,0)}(p_0) \neq 0$.
Since $\a^{(2,0)}(\d_z,\d_z)=\a_g(\d_z, \d_z)$ is isotropic, we have that
$\a^{*(2,0)}(\d_z,\d_z)$ is also isotropic. Define an $N_g L$-valued
tensor field on $U$ by
$$
\a^*=\a^{*(2,0)}+\overline{\a^{*(2,0)}}.
$$
By definition, the (1,1)-part of $\a^*$ vanishes, hence  it maps the unit
tangent circle at each tangent plane into an ellipse which, in fact, is a
circle of positive radius  since $\a^{*(2,0)}(\d_z,\d_z)$ is isotropic.
Now we may extend $N_1^g|_{L_0}$ to a plane bundle  $N_1^g$ defined over 
all $L^2$ by defining 
$$
N_1^g(p_0)=\text{span}\{\text{image}\,\a^*(p_0)\},
$$ 
and this concludes the proof.\vspace{1,5ex}\qed

To conclude this section, we show how to construct any $1$-isotropic  
simply connected surface in parametric form. This procedure can
easily be used to construct complete examples as was done in a quite similar 
situation in \cite{dg3}.\vspace{1ex}

On a simply connected domain $U\subset\C$, a minimal surface 
$g\colon U\to\R^N$  has the generalized Weierstrass representation
$$
g=\mbox{Re}\int^{\displaystyle{z}}\gamma dz
$$
where the Gauss map $\gamma\colon  U\to\C^N$ of $g$  has the expression
$$
\gamma=\frac{\beta}{2}\left(1-\phi^2,i(1+\phi^2),2\phi\right)
$$
being $\beta$ holomorphic and $\phi\colon U\to\C^{N-2}$ meromorphic;
see \cite{ho} for details. From \cite{Ch} we have that $g$ is $1$-isotropic  
if and only if $(\phi',\phi')=0$,
where $(\,,\,)$ stands for the standard symmetric inner product in $\C^{N-2}$.
Hence, to construct any $1$-isotropic  surface start with a nonzero holomorphic 
map $\alpha_0\colon  U\to\C^{N-4}$. 
Assuming that $\alpha_1\colon  U\to\C^{N-2}$ has been defined already, set
$$
\alpha_2 = \beta_2\left(1-\phi_1^2,i(1+\phi_1^2),2\phi_1\right)
$$
where $\phi_1=\int^z\alpha_1dz$ and $\beta_2\neq 0$ is any holomorphic function. 
Then, the  surface with Gauss map $\gamma=\alpha_1$, i.e., 
$g=\mbox{Re}\;\alpha_2$, is $1$-isotropic.

\section{The results}

In this section, we state the results of this paper and leave the proofs for
the following one. \vspace{1,5ex}

Let $g\colon L^2\to\R^{n+2}$, $n\geq 3$, be a substantial $1$-isotropic surface 
and let $\pi\colon\Lambda_g\to L^2$ denote the vector bundle of rank $n-2$ whose 
fibers are the orthogonal complement in the normal bundle $N_gL$ of $g$ of the extended
first normal bundle $N_1^g$ of $g$.  Associated to $g$ we consider the immersion 
$F_g\colon\Lambda_g\to\R^{n+2}$ given by 
\be\label{param}
F_g(p,v)=g(p)+v,
\ee
and denote by $M^n$ the manifold $\Lambda_g$ when it is endowed with the metric 
induced by $F_g$. By construction $F_g\colon M^n\to\R^{n+2}$ is an $(n-2)$-ruled
submanifold with complete rulings, that is, there is an integrable tangent distribution 
of dimension $n-2$ whose leaves are mapped diffeomorphically by $F$ onto complete 
affine subspaces of the ambient space.
\vspace{1ex}

In the sequel, we denote by ${\cal H}$ the tangent distribution 
orthogonal to the rulings. An embedded surface $j\colon L^2\to M^n$ 
is called an \emph{integral surface} of ${\cal H}$ if $j_*T_pL={\cal H}(j(p))$ 
at every point $p\in L^2$.

\begin{theorem}\po\label{main1} Let $g\colon L^2\to\R^{n+2}$, $n\geq 4$, be
a  $1$-isotropic substantial surface. Then the associated  immersion 
$F_g\colon M^n\to\R^{n+2}$ is an $(n-2)$-ruled minimal submanifold with rank  
$\rho=4$ on an open dense subset of $M^n$. Moreover, the rulings of $F_g$ are 
complete and the integral surface $L^2$ of ${\cal H}$ is unique and totally 
geodesic. Furthermore, the metric of $M^n$ is complete if $L^2$ is complete.

Conversely, let  $F\colon M^n\to\R^{n+2}$, $n\geq 4$, be an $(n-2)$-ruled 
minimal immersion with rank $\rho=4$ on an open dense subset of $M^n$.
Assume that ${\cal H}$ admits a  totally geodesic integral surface 
$j\colon L^2\to M^n$ which is a global cross section to the rulings. 
Then, the surface $g=F\circ j\colon L^2\to\R^{n+2}$ is $1$-isotropic and 
$F$ can be parametrized by (\ref{param}).
\end{theorem}
 
The vertical bundle $\mathcal V=\ker\pi_*$ of the submersion $\pi$ decomposes 
orthogonally as 
$$
\mathcal{V}=\mathcal{V}^1\oplus \mathcal{V}^0
$$ 
on an open dense subset of $L^2$, where $\mathcal{V}^1$ denotes the plane 
bundle determined by $N_2^g$. We assume without loss of
generality that this decomposition holds globally.
In the sequel, we consider the orthogonal decomposition of the tangent bundle
of $M^n$ given by $TM={\cal H}\oplus\mathcal{V}$
where we identify isometrically (and use the same
notation) the subbundle $\mathcal{V}$ tangent to the rulings with the 
corresponding normal subbundle to $g$. Then, it follows from the proof
that the relative nullity leaves of $F$ are identified with the fibers 
of $\mathcal{V}^0$.
\medskip

Let $\cal J$ be the endomorphism of $TM$ such that 
${\cal J|_{\cal{H}}}\colon\cal H\to\cal H$ is the almost complex 
structure in $\mathcal{H}$ determined by the orientation and restricted 
to $\cal V$ is the identity, and set
$$
{\cal J}_{\theta}=\cos\theta I+\sin\theta {\cal J}.
$$

\begin{theorem}\po\label{main2}  Let $g\colon L^2\to\R^{n+2}, n\geq 4$, 
be a simply connected $1$-isotropic substantial surface. Then $F_g\colon M^n\to\R^{n+2}$ 
allows a smooth one-parameter family of minimal \mbox{genuine} isometric deformations
$F_\theta\colon M^n\to\R^{n+2},\;\theta\in\Sf^1$, 
such that $F_0=F_g$ and each $F_\theta$ carries the same ruling and 
relative nullity leaves as $F_g$.  

Moreover, there is  a parallel vector bundle isometry 
$\Psi_\theta\colon N_{F_g}M\to N_{F_\theta}M$
such that the relation between the second 
fundamental forms of $F_\theta$ and $F_g$ is given by
\be\label{forms}
\a_{F_\theta}(X,Y)=\Psi_\theta \big(\textsf{R}_{-\theta} \a_{F_g}(X,Y) 
+2\kappa\sin(\theta/2)\beta ({\cal J}_{-\theta /2}X,Y)\big)
\ee
where $\textsf{R}_\theta$ is the rotation of angle
$\theta$ on $N_{F_g}M$ that preserves orientation,
$\kappa$ is the radius of the ellipse of curvature of $g$
and  $\beta$ is the traceless bilinear form  defined by (\ref{b}).
\end{theorem}

\begin{remark}\po {\em Quite similar arguments give that the above two 
results hold for dimension $n=3$ and rank $\rho=3$.
}\end{remark}

If $g$ is holomorphic with respect to some parallel complex structure 
in $\R^{n+2}$, then taking  a rotation of angle $\theta$ that  
preserves orientation in each $N_s^g$, $s\geq 2$, induces 
an intrinsic isometry $S_\theta$ on $ M^n$.

\begin{theorem}\po\label{main3}  If $g$ is holomorphic then $F_g\circ S_{-\theta}$ 
is congruent to $F_{\theta}$ for any $\theta\in\Sf^1$.  
\end{theorem}

\section{The proofs}

Let  $g\colon L^2\to\R^{n+2},n\geq 4,$ be a substantial oriented minimal 
surface. We  choose local positively oriented orthonormal frames 
$\{e_1,e_2\}$ in $TL$ and $\{e_3,e_4\}$ in $N_1^g$ 
such that 
$$
\a_g(e_1,e_1)=\kappa e_3\;\;\;\mbox{and}\;\;\; \a_g(e_1,e_2)=\mu e_4
$$ 
where $\kappa,\mu$ are the semi-axes of the ellipse of curvature. We also take 
a local orthonormal normal frame $\{e_5,\ldots,e_{n+2}\}$ 
such that $\{e_{2r+1},e_{2r+2} \}$ is a positively oriented frame field spanning 
$N_r ^g$ for every even $r$. When $n=2m+1$ is odd, then $e_{2m+1}$ spans the last 
normal bundle. We refer to  $\{e_1,\ldots,e_{n+2}\}$ as an \emph{adapted frame} 
of $g$ and consider the one-forms
$$
\omega_{ij}=\<\nab e_i,e_j\>\;\;\mbox{for}\;\;1\leq i,j\leq n+2.
$$
Then, we have from
$$
\a_g^3(e_1,e_1,e_1)+\a_g^3(e_1,e_2,e_2)=0\;\;\;\mbox{and}\;\;\;
\a_g^3(e_1,e_1,e_2)=\a_g^3(e_2,e_1,e_1)
$$
that
\be\label{conn}
\omega_{45}=-\dfrac{1}{\lambda}*\omega_{35}
\;\;\mbox{and}\;\;\omega_{46}=-\dfrac{1}{\lambda}*\omega_{36}
\ee
where  $\lambda= \mu/\kappa$,  $*$ denotes the Hodge operator, i.e., 
$*\omega(e)=-\omega(Je)$, and  $J$ is the  complex structure of $L^2$ induced
by the metric and the orientation.  We denote by 
$$
V=a_1e_1+a_2e_2,\; W=b_1e_1+b_2e_2,\;
Y=c_1e_1+c_2e_2\;\;\mbox{and}\;\;
Z=d_1e_1+d_2e_2
$$
the dual vector fields of $\omega_{35},\omega_{36},\omega_{45}$ and $\omega_{46}$, 
respectively.  Then (\ref{conn}) is equivalent to 
$$
Y=-\frac{1}{\lambda}JV \;\;\mbox{and}\;\; Z=-\frac{1}{\lambda}JW,
$$
and hence
$$
\lambda c_1 = a_2,\; \lambda c_2=-a_1,\; 
\lambda d_1 = b_2\;\;\mbox{and}\;\;\lambda d_2=-b_1.
$$
Clearly $F=F_g$ is an immersion and the horizontal bundle $\mathcal H$ is the  
orthogonal complement of $\mathcal V$ in the tangent bundle of $M ^n$, i.e., 
we have at $(p,v)\in M ^n$ that
$$
T_{(p,v)}M=\mathcal H(p,v)\oplus\mathcal V(p,v).
$$ 

Fixed $(p,v)\in M ^n$, define a normal vector field $\delta_v$ in a 
neighborhood $U$ of $p$ by 
\be\label{delta}
\delta_v(q)=\sum_{j\geq 5}\<v,e_j(p)\> e_j(q).
\ee  
Let $\beta_i$, $1\leq i\leq 2$, be the curves in $M^n$ 
satisfying $\beta_i(0)=(p,v)$  given by 
$$
\beta_i(s)=(c_i(s), \delta_v(c_i(s)))
$$
where each $c_i(s)$ is a smooth curve in $L^2$ such that $c_i'(0)=e_i(p)$. 
Then $Y_1,Y_2\in T_{(p,v)}M$ where
\be\label{yes}
Y_i= \beta_i'(0),\;1\leq i\leq 2.
\ee
Let $G_i,H_i\in C^\infty(M)$, $\;1\leq i\leq 2$, be the functions 
$$
G_i=t_2\omega^i_{56}+t_3\omega^i_{57}+t_4\omega^i_{58},\;\;
H_i=-t_1\omega^i_{56}+t_3\omega^i_{67}+t_4\omega^i_{68}
$$
where $\omega_{ij}^k=\omega_{ij}(e_k)$ and 
$t_j\in C^\infty(M)$ is defined by
$$
t_j(q,w)=\<w,e_{j+4}(q)\>,\;\; 1\leq j\leq 4.
$$

The vertical bundle $\mathcal V$   can be orthogonally decomposed as
$\mathcal{V}=\mathcal{V}^1\oplus \mathcal{V}^0$
where $\mathcal{V}^1$ denotes the plane bundle determined by $N_2^g$. 
Let $\{E_3, E_4\}$ and $\{E_5,\dots, E_n\}$ be the local orthonormal frames  
of $\mathcal V^1$ and  $\mathcal V^0$, respectively, such that 
$$
F_* E_j= e_{j+2},\;\;3\leq j\leq n.
$$

\begin{lemma}\label{NF}\po The vectors  $X_1,X_2\in T_{(p,v)}M$ 
defined as  
\be\label{X}
X_i= Y_i+ G_iE_3+H_iE_4-\sum_{j\geq 7}\<\nap_{e_i}\delta_v,e_j\>E_{j-2}
\ee
satisfy $X_1,X_2\in\mathcal H(p,v)$ and 
$$
F_*X_1=g_* e_1 -\va_1e_3
-\dfrac{1}{\lambda}\va_2e_4,\;\;\;
F_*X_2=g_* e_2 -\va_2e_3
+\dfrac{1}{\lambda}\va_1e_4
$$
where 
$$
\va_j=t_1^0a_j+t_2^0b_j,\; j=1,2,
$$
and $t_j^0=t_j(p,v)$.
Moreover, the normal space  $N_FM(p,v)$ is spanned by 
$$
\xi=g_*(t_1^0V(p)+ t_2^0W(p))+e_3(p),\;\;
\eta=g_*(t_1^0Y(p)+ t_2^0Z(p))+e_4(p).
$$
In particular, if $g$ is $1$-isotropic we have
$$
\|X_1\|^2 = \|X_2\|^2=\Omega^2=1+\|t_1^0V(p)+t_2^0W(p)\|^2 
,\;\; \<X_1,X_2\>=0
$$
and 
$$
\|\xi\|=\Omega=\|\eta\|,\;\; \<\xi,\eta\>=0.
$$
\end{lemma}

\proof We obtain from 
$$
F_* Y_i=g_{*_p}e_i(p)
+\sum_{j\geq3}\<\nap_{e_i}\delta_v,e_j\>(p)e_j(p)
$$
that
$$
F_*Y_i-\sum_{j\geq5}\<\nap_{e_i}\delta_v, e_j\>(p)F_*E_{j-2}
=g_{*_p}e_i(p)-\sum_{3\leq k\leq 4}\<\nap_{e_i}e_k,\delta_v\>(p)e_k(p).
$$
On the other hand, 
$$
\<\nap_{e_i}\delta_v,e_5\>(p)= -t_2^0\omega^i_{56}(p)
-t_3^0\omega^i_{57}(p)-t_4^0\omega^i_{58}(p)=-G_i(p,v),
$$
$$
\<\nap_{e_i}\delta_v, e_6\>(p)= t_1^0\omega^i_{56}(p)
-t_3^0\omega^i_{67}(p)-t_4^0\omega^i_{68}(p)=-H_i(p,v)
$$
and
$$
\<\nap_{e_i}e_3,\delta_v\>(p)=t_1^0\omega^i_{35}(p)+t_2^0\omega^i_{36}(p)
=t_1^0a_i(p)+t_2^0b_i(p),
$$
$$
\<\nap_{e_i}e_4,\delta_v\>(p)=t_1^0\omega^i_{45}(p)+t_2^0\omega^i_{46}(p)
=t_1^0c_i(p)+t_2^0d_i(p)
$$
where also $t_j^0=t_j(p,v),\;3\leq j\leq 4$. Hence, 
$$
F_*X_i=g_* e_i -(t_1^0a_i+t_2^0b_i)e_3
-(t_1^0c_i+t_2^0d_i)e_4,\;i=1,2.
$$
The remaining of the proof is immediate.
\qed

\begin{lemma}\po\label{comp}
The following equations hold:
\be\label{first}
\xi_* E_3=g_*V,\;\;\xi_* E_4=g_*W\;\;\mbox{and}\;\;\xi_*=0\;
\text{on}\;\mathcal{V}^0,
\ee
\be\label{second}
\eta_*E_3=g_*Y,\;\;\eta_* E_4=g_*Z\;\;\mbox{and}\;\;
\eta_*=0\;\text{on}\;\mathcal{V}^0,
\ee
\begin{eqnarray}\label{1}
\xi_* X_1\!\!\!&=&\!\!\!g_*\big((e_1(\va_1)-\kappa)e_1+e_1(\va_2)e_2
+\omega_{12}^1J(t_1V+t_2W)+G_1V+H_1W\big)\nonumber\\
\!\!\!&&\!\!\!+\kappa\va_1e_3+(\omega_{34}^1+\lambda\kappa\va_2)e_4+a_1e_5+b_1e_6,
\end{eqnarray}
\begin{eqnarray}
\xi_* X_2\!\!\!&=&\!\!\!g_*\big(e_2(\va_1)e_1+(e_2(\va_2)+\kappa)e_2
+\omega_{12}^2J(t_1V+t_2W)+G_2V+H_2W\big)\nonumber\\
\!\!\!&&\!\!\!-\kappa\va_2e_3+(\omega_{34}^2+\lambda\kappa\va_1)e_4+a_2e_5+b_2e_6,
\end{eqnarray}
\begin{eqnarray}
\eta_* X_1\!\!\!&=&\!\!\!g_*\big(e_1(\psi_1)e_1+(e_1(\psi_2)-\lambda\kappa)e_2
+\sigma\omega_{12}^1(t_1V+t_2W)-\sigma G_1JV-\sigma H_1JW\big)\nonumber\\
\!\!\!&&\!\!\!-(\omega_{34}^1-\kappa\psi_1)e_3
+\lambda\kappa\psi_2e_4+\sigma a_2e_5+\sigma b_2e_6,
\end{eqnarray}
\begin{eqnarray}\label{4}
\eta_* X_2\!\!\!&=&\!\!\!g_*\big((e_2(\psi_1)-\lambda\kappa)e_1+e_2(\psi_2)e_2
+\sigma\omega_{12}^2(t_1V+t_2W)-\sigma G_2JV-\sigma H_2JW
\big)\nonumber\\
\!\!\!&&\!\!\!-(\omega_{34}^2+\kappa\psi_2)e_3
+\lambda\kappa\psi_1e_4-\sigma a_1e_5-\sigma b_1e_6
\end{eqnarray}
where $\sigma=1/\lambda$ and 
$\psi_j=t_1^0c_j+t_2^0d_j,\, j=1,2$.
\end{lemma}

\proof Let $\gamma(s)=(c(s), v(s))$ be a curve in $M^n$ so that 
$\gamma(0)=(p,v)$ and $\gamma'(0)\in\mathcal{V}(p,v)$, 
that is, $c'(0)=0$. We have that
$$
\xi_*\gamma'(0)=\<Dv/ds(0),e_5(p)\>g_*V(p)+\<Dv/ds(0),e_6(p)\> g_*W(p), 
$$
or equivalently, that
$$
\xi_*\gamma'(0)=\<F_*\gamma'(0),e_5(p)\>g_*V(p)
+\<F_*\gamma'(0),e_6(p)\> g_*W(p). 
$$
 From this we obtain (\ref{first}). Similarly, we have
$$
\eta_*\gamma'(0)=\<F_*\gamma'(0),e_5(p)\>g_*Y(p)
+\<F_*\gamma'(0),e_6(p)\>g_*Z(p)
$$
from which we obtain (\ref{second}). 

Making use of Lemma \ref{NF} and the Gauss and Weingarten  formulas for $g$
we compute equations (\ref{1}) to (\ref{4}). 
We only argue for  (\ref{1})  since the proof of the other  equations 
is completely  similar.  We have from (\ref{X}) and (\ref{first}) that
$$
\xi_*X_i =\xi_*Y_i+G_ig_*V +H_ig_*W,\;\;1\leq i\leq 2. 
$$
In view of (\ref{yes}) and since
$$
(\xi\circ \beta_i)(s)=t_1^0g_*V(c_i(s))
+t_2^0g_*W(c_i(s))+e_3(c_i(s))
$$
we obtain
\bea
\xi_*Y_1\!\!\!&=&\!\!\!t_1^0\big(g_*\n_{e_1}V +\a_g(e_1,V)\big)(p)
+t_2^0\big(g_*\n_{e_1}W +\a_g(e_1,W)\big)(p)\\
\!\!\!&&\!\!\!-\kappa(p)g_*e_1(p)+\nap_{e_1}e_3(p),
\eea
and the desired formula for $\xi_*X_1$ follows by direct computations.
\qed

\begin{lemma}\po\label{A}
If $g$ is a $1$-isotropic surface, then the shape operators of $F_g$ with 
respect to the orthonormal tangent frame 
$$
E_i=X_i/\Omega,\;i=1,2,\;\;\text{and}\;\;F_*E_j=e_{j+2},
\;3\leq j\leq n,
$$
vanish along $\mathcal{V}^0$ 
and restricted to $\mathcal{H}\oplus\mathcal{V}^1$  are given by
\be\label{ssf}
A_{\xi}=\begin{bmatrix}
\kappa + h_1&h_2&r_1&s_1\\
h_2&-\kappa-h_1&r_2&s_2\\
r_1&r_2&0&0\\
s_1&s_2&0&0&
\!\!\!\!\!\end{bmatrix},\;\;\;
A_{\eta}=\begin{bmatrix}
h_2&\kappa-h_1&r_2&s_2\\
\kappa-h_1&-h_2&-r_1&-s_1\\
r_2&-r_1&0&0\\
s_2&-s_1&0&0&
\!\!\!\!\!\end{bmatrix}\\
\ee
where $r_i\Omega=-a_i$, $s_i\Omega=-b_i$, 
\bea
h_i\!\!\!&=&\!\!\!-\dfrac{1}{\Omega^2}\big(t_1(e_i(a_1)
- a_2B_i- b_1\omega_{56}^i) 
+ t_2(e_i(b_1)- b_2B_i+ a_1\omega_{56}^i)\\
\!\!\!&&\!\!\!+\,t_3(a_1\omega_{57}^i+ b_1\omega_{67}^i)
+t_4(a_1\omega_{58}^i+b_1\omega_{68}^i)\big)
\eea
and $B_i=\omega_{12}^i+\omega_{34}^i$, $i=1,2$.
\end{lemma}

\proof Since $g$ is $1$-isotropic, then (\ref{1}) to (\ref{4}) 
hold for $\psi_1=\va_2$ and $\psi_2=-\va_1$.
On the other hand, a straightforward computation shows that 
the Ricci equations 
$$
\<R^\perp(e_1,e_2) e_\a,e_\beta\>=0
$$ 
for $\a=3,4$ and $\beta=5,6$ are equivalent to
\bea
&&e_1(a_2)-e_2(a_1)+a_1B_1+a_2B_2-b_2\omega_{56}^1 
+b_1\omega_{56}^2=0,\\
&&e_1(b_2)-e_2(b_1)+b_1B_1+b_2B_2+a_2\omega_{56}^1
-a_1\omega_{56}^2=0,\\
&&e_1(a_1)+e_2(a_2)-a_2B_1+a_1B_2-b_1\omega_{56}^1
-b_2\omega_{56}^2=0,\\
&&e_1(b_1)+e_2(b_2)-b_2B_1+b_1B_2+a_1\omega_{56}^1
+a_2\omega_{56}^2=0,
\eea
and for $\a=3,4$ and $\beta=7,8$ are equivalent to
\bea
&&a_2\omega_{57}^1-a_1\omega_{57}^2 
+b_2\omega_{67}^1-b_1\omega_{67}^2=0,\\
&&a_2\omega_{58}^1-a_1\omega_{58}^2 
+b_2\omega_{68}^1-b_1\omega_{68}^2=0,\\
&&a_1\omega_{57}^1+a_2\omega_{57}^2 
+b_1\omega_{67}^1+b_2\omega_{67}^2=0,\\
&&a_1\omega_{58}^1+a_2\omega_{58}^2 
+b_1\omega_{68}^1+b_2\omega_{68}^2=0.
\eea
We thus have that
$$
\<A_\xi E_i,E_j\>=-\<F_*E_i,\xi_* E_j \>\;\;\mbox{and}\;\; 
\< A_\eta E_i,{E_j}\>=-\< F_*E_i ,\eta_* E_j \>,\;1\leq i,j \leq n,
$$ 
and the result follows by a straightforward computation.\vspace{1,5ex}\qed

\noindent\emph{Proof of Theorem \ref{main1}:}
We first prove the converse. Let $F\colon M^n\to\R^{n+2}$, $n\geq 4$, 
be an $(n-2)$-ruled minimal immersion with rank $\rho=4$ on 
an open dense subset. Then the tangent bundle splits as $TM={\cal H} \oplus {\cal V}$, 
where $\cal H$ is orthogonal to the rulings and ${\cal V}$ splits as
${\cal V}={\cal V}^1\oplus{\cal V}^0$ with the fibers of ${\cal V}^0$ being 
the relative nullity leaves.  

 The normal space of the surface  $g=F\circ j$  at any point $x\in L^2$ is 
given by 
$$
N_gL(x)=F_*(j(x)){\cal V}\oplus N_FM(j(x)).
$$
Let $\Lambda_g$ be the subbundle of the normal bundle of $g$ whose fiber 
at $x\in L^2$ is $F_*(j(x)){\cal V}$. 
Observe that 
$$
F(p)-g\circ\pi (p)=F(p)-F(j(x))\in F_*(j(x)){\cal V}
$$ 
for any $p\in M^n$, where $x=\pi(p)$, since $p$ and $j(x)$ belong to the 
same leaf of  ${\cal V}$. 
Since $F$ maps diffeomorphically the leaves of ${\cal V}$  
onto complete  affine subspaces, it follows that the map 
$T\colon M^n\to\Lambda_g$ given by
$$
T(p)=(\pi(p),F(p)-g\circ\pi(p)) 
$$
is  a global diffeomorphism.  Clearly the  immersion 
$\tilde{F}=F\circ T^{-1}$ satisfies 
$$
\tilde{F}(x,v)=g(x)+v,
$$
i.e., $\tilde{F}=F_g$ is of the form  (\ref{param}). 
Identifying $M^n$ with $\Lambda_g$ via $T$, we have 
that $F=F_g$ and $j$ is the zero section of $\Lambda_g$.

It remains to show that $g$ is 1-isotropic. 
Being $j$  totally geodesic, we have that
\be\label{ag}
\a_g(X,Y)=\a_F (j_*X, j_*Y)
\ee
for all $X,Y\in TL$. This and our assumptions imply that $g$ is minimal. 
The horizontal and the vertical bundles satisfy
$$
F_*(p,v)\mathcal V=(N_1^g(p))^\perp\subset N_gL(p),\;\;
F_*(p,v)\mathcal H\subseteq g_*T_pL\oplus(\Lambda_g(p))^\perp,
$$
$$
N_FM(p,v)\subseteq g_*T_pL\oplus(\Lambda_g(p))^\perp
$$
and now (\ref{ag}) yields  $N_1^g=\Lambda_g^\perp$.   

Let $\{e_1,\ldots,e_{n+2}\}$ be an adapted frame of $g$.
Setting 
$$
g_{ij}=\<X_i,X_j\>_F,\; b^\xi_{ij}=\<\xi_*X_i,F_*X_j\>\;\;\mbox{and}\;\; 
b^\eta_{ij}=\<\eta_*X_i,F_*X_j\>,\; i,j=1,2,
$$
and using Lemma \ref{NF} and Lemma \ref{comp}, we find that
$$
g_{11}=1+\va_1^2+\sigma^2\va_2^2,\;\;
g_{12}=(1-\sigma^2)\va_1\va_2,\;\;
g_{22}=1+\va_2^2+\sigma^2\va_1^2,
$$
and
\bea
b^\xi_{11}\!\!\!&=&\!\!\!e_1(\va_1)-\kappa-\omega_{12}^1\va_2+G_1a_1+H_1b_1
-\kappa\va_1^2-\sigma\va_2(\omega_{34}^1+\mu\va_2),\\
b^\xi_{12}\!\!\!&=&\!\!\!e_1(\va_2)+\omega_{12}^1\va_1+G_1a_2+H_1b_2
-\kappa\va_1\va_2 +\sigma\va_1(\omega_{34}^1+\mu\va_2),\\
b^\xi_{21}\!\!\!&=&\!\!\!e_2(\va_1)-\omega_{12}^2\va_2+G_2a_1+H_2b_1
+\kappa\va_1\va_2-\sigma\va_2(\omega_{34}^2+\mu\va_1),\\
b^\xi_{22}\!\!\!&=&\!\!\!e_2(\va_2)+\kappa+\omega_{12}^2\va_1+G_2a_2+H_2b_2
+\kappa\va_2^2+\sigma\va_1(\omega_{34}^2+\mu\va_1)
\eea
and
\bea
b^\eta_{11}\!\!\!&=&\!\!\!e_1(\psi_1)-\omega_{12}^1\psi_2+\sigma G_1a_2+\sigma H_1 b_2
+\omega_{34}^1\va_1-\kappa(\va_1\psi_1+\va_2\psi_2),\\
b^\eta_{12}\!\!\!&=&\!\!\!e_1(\psi_2)-\mu+\omega_{12}^1\psi_1- \sigma G_1a_1-\sigma H_1b_1
+\omega_{34}^1\va_2+\kappa(\va_1\psi_2-\va_2\psi_1),\\
b^\eta_{21}\!\!\!&=&\!\!\!e_2(\psi_1)-\mu-\omega_{12}^2\psi_2+\sigma G_2a_2+\sigma H_2b_2
+\omega_{34}^2\va_1+\kappa(\va_1 \psi_2-\va_2\psi_1),\\
b^\eta_{22}\!\!\!&=&\!\!\!e_2(\psi_2)+\omega_{12}^2\psi_1-\sigma G_2a_1-\sigma H_2b_1
+\omega_{34}^2\va_2 +\kappa(\va_1\psi_1+\va_2\psi_2).
\eea
 From our assumptions, we have 
\be\label{m1}
g_{11}b^\xi_{22}-g_{12}(b^\xi_{12}+b^\xi_{21})+g_{22}b^\xi_{11}=0
\ee
and
\be\label{m2}
g_{11}b^\eta_{22}-g_{12}(b^\eta_{12}+b^\eta_{21})+g_{22}b^\eta_{11}=0.
\ee
Viewing (\ref{m1}) and (\ref{m2})  as polynomials were the coefficients of 
$t_1^4, t_2^4, t_1^2t_2^2$  must vanish gives
$$
(\lambda^2-1)(a_1^2+a_2^2)(a_1^2-a_2^2)=0=(\lambda^2-1)(b_1^2+b_2^2)(b_1^2-b_2^2)
$$
and
$$
(\lambda^2-1)a_1a_2(a_1^2+a_2^2)=0=(\lambda^2-1)b_1b_2(b_1^2+b_2^2).
$$
Hence $\lambda=1$ since otherwise, we have from the above that
$\omega_{35}=\omega_{36}=\omega_{45}=\omega_{46}=0$, 
which is a contradiction.\vspace{1ex} 

We now prove the direct statement. Observe that $g=F_g\circ j$, where $j$ is the 
zero section of $M^n$.  Clearly, we have that $j$  is an  integral surface of the 
distribution orthogonal to the rulings which is also totally geodesic and a global 
cross section to the rulings. Up to the uniqueness of the integral surface and 
completeness of $M^n$ the  proof now follows from Lemma \ref{A}. In fact,
it is very easy to see that the metric of $M^n$ is complete if the metric 
of $L^2$ is complete. 

Assume  that there exists a second integral surface $\tilde j\colon L^2\to M^n$. 
Set $\tilde g=F_g\circ\tilde j$ and let  $\tilde T\colon M^n\to\Lambda_{\tilde g}$ 
be the diffeomorphism given by 
$$
\tilde T(p)=(\pi(p),F(p)-\tilde g(\pi(p)).
$$
Then $\tilde T\circ T^{-1}\colon\Lambda_g\to\Lambda_{\tilde g}$ is 
$$
\tilde T \circ T^{-1} (x,v)=(x,v+g(x)- \tilde g(x)).
$$
Hence $\Lambda_g$ and $\Lambda_{\tilde g}$  can be identified by parallel translation, 
thus there exists a section $\delta$ of $\Lambda_g$ such that $\tilde g=g+\delta$. 
It follows from
\be\label{d}
\tilde g_*X=g_*X+\nap_X \delta
\ee
that $\nap_X\delta\in N_1^g$ for any $X\in TL$. If $\delta$ is 
constant, then $g$ lies  in an affine subspace $\R^{n+1}$ of $\R^{n+2}$ 
perpendicular to $\delta$ which has been excluded.  Thus, there is 
$\mu=\nap_{X_0}\delta\neq 0$ for some $X_0\in TL$. From (\ref{d}) we have  that
$\nap_Y\mu\in N_1^g$ for any $Y\in TL$. This easily implies that $ N_1^g$ 
is parallel in the normal bundle and thus $g$ lies in $\R^4$, a contradiction.
\vspace{1,5ex}\qed

\noindent\emph{Proof of Theorem \ref{main2}:} 
For each $\theta\in\Sf^1$, we define $F_\theta\colon\Lambda_g\to\R^{n+2}$ by
$$
F_\theta(p,v)=g_\theta(p)+v.
$$ 
In the sequel, corresponding quantities of $F_\theta$ are denoted 
by the same symbol used for $F_g$ marked with $\theta$.  
That $F_\theta$ is isometric to $F_g$ is immediate.
Since the tangent frame $\{e_1,e_2\}$
has been fixed, we have for the adapted frames of $g_\theta$ that
$$
e^\theta_3=R^1_\theta e_3\;\;\mbox{and}\;\; e^\theta_4=R^1_\theta e_4        
$$
where $R^1_\theta$ is the  rotation of angle $\theta$ on $N^g_1$. 
We complete the adapted frame choosing 
$$
e^\theta_j=e_j,\;\;\; 5\leq j \leq n+2.
$$
Clearly, it holds that $\omega^\theta_{34}=\omega_{34}$ 
and $\omega^\theta_{ij}=\omega_{ij}$ for $i,j\geq 5$. Moreover, 
$$
\omega^\theta_{35}
=\cos\theta\omega_{35}-\sin\theta*\omega_{35}\;\;\mbox{and}\;\;
\omega^\theta_{36}=\cos\theta\omega_{36}-\sin\theta*\omega_{36}.
$$
Hence, the dual vector fields of $\omega^\theta_{36}$  
and $\omega^\theta_{36}$ are given, respectively, by
$$
V_\theta =J_{-\theta}V \;\; \text{and}\;\; W_\theta =J_{-\theta}W.
$$
Thus,
$$
a_1^\theta= a_1\cos\theta+  a_2\sin\theta,\;\;a_2^\theta
= a_2\cos\theta- a_1\sin\theta
$$
and 
$$
  b_1^\theta=  b_1\cos\theta+  b_2\sin\theta, \;\;   b_2^\theta
=  b_2\cos\theta-  b_1\sin\theta.
$$
It follows from (\ref{delta}), (\ref{yes}) and  (\ref{X}) that 
$$
X_i^\theta=X_i,\;\;i=1,2.
$$
By Lemma \ref{NF}, the normal bundle of $F_\theta$ is spanned by 
$$
\xi_\theta=g_{\theta_*}J_{-\theta}(t_1V+ t_2W)+R^1_\theta e_3,\;\;\;
\eta_\theta=-g_{\theta_*}J_{\pi/2-\theta}(t_1V+t_2W)+R^1_\theta e_4.
$$
A straightforward computation yields that the map
$\Psi_\theta\colon N_{F_g}M\to N_{F_\theta}M$ given by
$$
\Psi_\theta \xi=\xi_\theta \;\;\text{and}\;\;\Psi_\theta\eta=\eta_\theta
$$
is a parallel vector bundle isometry.  The shape operators of 
$F_\theta$ vanish on $\mathcal{V}^0$
and restricted to $\mathcal{H}\oplus\mathcal{V}^1$ are given with respect to the 
frame $\{E_1,\dots,E_n\}$  by
$$
A^\theta_{\xi_\theta}=\begin{bmatrix}
\kappa + h_1^\theta&h_2^\theta&r_1^\theta&s_1^\theta\\
h_2^\theta&-\kappa-h_1^\theta&r_2^\theta&s_2^\theta\\
r_1^\theta&r_2^\theta&0&0\\
s_1^\theta&s_2^\theta&0&0&
\!\!\!\!\!\end{bmatrix}, 
\;\;\;
A^\theta_{\eta_\theta}=\begin{bmatrix}
h_2^\theta&\kappa-h_1^\theta&r_2^\theta&s_2^\theta\\
\kappa-h_1^\theta&-h_2^\theta&-r_1^\theta&-s_1^\theta\\
r_2^\theta&-r_1^\theta&0&0\\
s_2^\theta&-s_1^\theta&0&0&
\!\!\!\!\!\end{bmatrix}
$$
where $r_i^\theta\Omega=-a_i^\theta$, 
$s_i^\theta\Omega=-b_i^\theta$ and 
$$
h_1^\theta=h_1 \cos\theta +h_2 \sin\theta,
\;\;\;h_2^\theta= -h_1\sin\theta+h_2\cos\theta.
$$
Let $L_\theta$ denote the endomorphism of $TM$  such that  $L_\theta |_{\cal{V}}=0$ and 
$L_\theta|_{\cal H}\colon\cal H\to\cal H$ is the reflection given 
by
$$
L_\theta|_{\cal H}=\begin{bmatrix}
-\sin(\theta/2)&\cos(\theta/2)\\
\cos(\theta/2)&\sin(\theta/2)&
\!\!\!\!\!\end{bmatrix}
$$
with respect to the tangent frame $\{E_1,E_2\}$. It follows easily  that
$$
A^\theta_{\Psi_\theta\xi}=A_{\textsf{R}_\theta\xi} 
-2\kappa\sin(\theta/2)L_\theta\;\;\mbox{and}
\;\;A^\theta_{\Psi_\theta\eta}=A_{\textsf{R}_\theta\eta}
-2\kappa\sin(\theta/2){\cal J}L_\theta.
$$
By a direct computation we obtain
$$
\a_{F_\theta}(X,Y)=\Psi_\theta \Big(\textsf{R}_{-\theta} \a_{F_g}(X,Y) 
-\frac{2\kappa}{\Omega^2}\sin(\theta/2) (\<L_\theta X,Y\>\xi
+ \<L_\theta {\cal J}X,Y\>\eta)\Big).
$$
Define  $\beta$ as the  symmetric section of $Hom(TM\times TM,N_{F_g}M)$
with nullity $\mathcal V$ such that
\be\label{b}
\beta(E_1,E_1)=\frac{1}{\Omega^2}\xi=-\beta(E_2,E_2)\;\;\mbox{and}\;\;\beta(E_1,E_2)
=-\frac{1}{\Omega^2}\eta,
\ee
and the proof of (\ref{forms}) follows easily.  

Finally, that the isometric deformation $F_\theta$ of $F_g$ is genuine 
is immediate from (\ref{ssf}) since the shape operators of $F_g$ 
have rank four for any normal direction along an open dense subset. 
\vspace{1,5ex}\qed

\noindent\emph{Proof of Theorem \ref{main3}:} 
Being $g$  holomorphic, there exists an isometry 
$\tau$ of $\R^{n+2}$  such that $g_\theta=\tau\circ g$. 
The higher fundamental forms satisfy 
$$
\a^{s+1}_{g_\theta}=\tau_*\circ\a^{s+1}_g\;\;\mbox{for any}\;\; s\geq 1.
$$
It was shown in \cite{df1} that the almost complex structure $J$  
induces an almost complex structure $J_s$ on each $N_s^g$ defined by
$$
J_s\alpha^{s+1}_g(X_1,\ldots,X_{s+1})=\alpha^{s+1}_g(JX_1,\ldots,X_{s+1}).
$$
In the present case each $J_s\colon N^g_s\to N^g_s$ is an isometry.
Thus, we have
$$
\a^{s+1}_{g_\theta}=R^s_{\theta}\circ\a^{s+1}_g,
$$
where $R^s_{\theta}=\cos\theta I+\sin\theta J_s$. 
Hence $R^s_{\theta}=\tau_*|_{N^g_s}$. 
It is now easy to see that
$
F_\theta=\tau\circ F_g\circ S_{-\theta},
$
and this concludes the proof.\qed

\section{The case of holomorphic curves}

Let the substantial surface $g\colon L^2\to\R^{n+2}$, $n\geq 6$, be a holomorphic 
curve with respect to some parallel complex structure in $\R^{n+2}$.  
Let $\{e_1,e_2\}$ be an orthonormal tangent frame  such that
$$
\a_g^{s+1}(e_1,\dots,e_1)=\kappa_se_{2s+1},\;\;\a_g^{s+1}(e_1,\dots,e_1,e_2)
=\kappa_se_{2s+2},\;\;1\leq s\leq n/2.
$$
Then, set $\tau_s=\kappa_s/\kappa_{s-1}$, $1\leq s\leq n/2$, with $\kappa_0=1$.
It is well-known that $\kappa_s$ can be defined  as the radius of the 
$s^{th}$-curvature ellipse (cf.\ \cite{dv}) and that the functions $\tau_s$ 
are completely determined by the metric of $L^2$ in an explicit form by a 
result of Calabi (cf.\ \cite{L}).
\vspace{1,5ex}

We see next that in this case of a holomorphic curve $g$ the second fundamental 
form of the associated minimal ruled submanifold $F_g\colon M^n\to\R^{n+2}$ is 
substantially simpler than in the general case and completely determined by 
the metric of the surface.

\begin{proposition}\po\label{Ap}
Let $g\colon L^2\to\R^{n+2}$, $n\geq 6$, be  holomorphic. Then the shape 
operators of $F_g\colon M^n\to\R^{n+2}$ with
respect to the orthonormal tangent frame
$$
E_i=X_i/\Omega,\;i=1,2,\;\;\text{and}\;\;F_*E_j=e_{j+2},\;\; 3\leq j\leq n.
$$
vanish along $\mathcal{V}^0$ and restricted to $\mathcal{H}\oplus\mathcal{V}^1$  
are given by

\be\label{ssfc}
A_{\xi}=\begin{bmatrix}
\tau_1 +h_1&h_2&r&0\\
h_2&-\tau_1-h_1&0&r\\
r&0&0&0\\
0&r&0&0&
\!\!\!\!\!\end{bmatrix},\;\;\;
A_{\eta}=\begin{bmatrix}
h_2&\tau_1-h_1&0&r\\
\tau_1-h_1&-h_2&-r&0\\
0&-r&0&0\\
r&0&0&0&
\!\!\!\!\!\end{bmatrix}\\
\ee
where
\bea
h_1\!\!\!&=&\!\!\!-\dfrac{1}{1+(t_1^2+t_2^2)\tau_2^2}\big(t_1e_1( \tau_2)
- t_2e_2( \tau_2)+t_3\tau_2\tau_3\big),\\
h_2\!\!\!&=&\!\!\!-\dfrac{1}{1+(t_1^2+t_2^2)\tau_2^2}\big(t_1e_2( \tau_2)
+t_2e_1( \tau_2)+t_4\tau_2\tau_3\big),\\
r\!\!\!&=&\!\!\!-\dfrac{\tau_2}{\sqrt{1+(t_1^2+t_2^2)\tau_2^2}}\cdot
\eea
Moreover, the second fundamental form of $F_g$ depends only 
on the metric of $L^2$.
\end{proposition}

\proof From the choice of the normal frame and  the definition of higher 
fundamental forms, we find that the normal connection forms
$$
\omega_{\a\beta}^j=\<\nap_{e_j} e_{\a},e_{\beta}\>,\;\;1\leq j\leq 2,\;\
3\leq \a,\beta\leq n+2,
$$
satisfy
\bea
\a^{s+1}_g(e_1,\dots,e_1) \!\!\!&=&\!\!\! 
(\nap_{e_1}\a^s_g(e_1,\dots,e_1))_{N_s^g}\\
\!\!\!&=&\!\!\!\kappa_{s-1}(\nap_{e_1}e_{2s-1})_{N_s^g}\\
\!\!\!&=&\!\!\!\kappa_{s-1}\big(\omega^1_{2s-1,2s+1}e_{2s+1}
+\omega_{2s-1,2s+2}^1e_{2s+2}\big).
\eea
Similarly, we find
\bea
\a^{s+1}_g(e_1,\dots,e_1, e_2) \!\!\!&=&\!\!\!
\kappa_{s-1}\big(\omega_{2s,2s+1}^1e_{2s+1}
+\omega_{2s,2s+2}^1e_{2s+2}\big),\\
\a^{s+1}_g(e_2,e_2,e_1\dots,e_1)\!\!\!&=&\!\!\!\kappa_{s-1}
\big(\omega_{2s,2s+1}^2e_{2s+1}+\omega_{2s,2s+2}^2e_{2s+2}\big),\\
\a^{s+1}_g(e_2,e_1\dots,e_1)\!\!\!&=&\!\!\!\kappa_{s-1}
\big(\omega_{2s-1,2s+1}^2e_{2s+1}+\omega_{2s-1,2s+2}^2e_{2s+2}\big).
\eea
Thus, we obtain
\be\label{con1}
\omega_{2s-1,2s+1}=\omega_{2s,2s+2}=\tau_s\omega_1,\;\;\omega_{2s-1,2s+2}
=-\omega_{2s,2s+1}=\tau_s\omega_2.
\ee
Moreover, from part $(ii)$ of Lemma $6$ in \cite{V} it follows that
\be\label{con2}
\omega_{2s+1,2s+2}=(s+1)\omega_{12}+*d\log\kappa_s,
\;\; 1\leq s \leq n/2.
\ee
Then, using (\ref{con1}), (\ref{con2}) we have from Lemma \ref{A} 
that the second fundamental form of $F_g$ is given by (\ref{ssfc}).\qed

\begin{remark}\po \em {Notice that in order to obtain  the expressions of the shape 
operators in the above  result we only used that the first three ellipses of
curvature are circles. In \cite{dv4} we will discuss when $M^n$ is Kaehler. 
}\end{remark}

\vspace{.5in} 
{\renewcommand{\baselinestretch}{1}
\hspace*{-20ex}\begin{tabbing} \indent\= IMPA -- Estrada Dona Castorina, 110
\indent\indent\= Univ. of Ioannina -- Math. Dept. \\
\> 22460-320 -- Rio de Janeiro -- Brazil  \>
45110 Ioannina -- Greece \\
\> E-mail: marcos@impa.br \> E-mail: tvlachos@uoi.gr
\end{tabbing}}
\end{document}